# Efficient calculation of transient eddy current response from multilayer cylindrical conductive media

Theodoros Theodoulidis

*Department of Mechanical Engineering*

*University of Western Macedonia*

*Bakola & Sialvera, 50132, Kozani*

*GREECE*

Email: theodoul@uowm.gr

Anastassios Skarlatos

*Département Imagerie Simulation pour le Contróle, Institute CEA LIST*

*CEA Saclay - Digiteo*

*Bat. 565-PC120*

*F-91191 Gif-sur-Yvette Cedex*

*FRANCE*

Email: Anastasios.SKARLATOS@cea.fr

**Abstract.** The transient response from a transmitter-receiver coil system inside a multilayer cylindrical conductive configuration is obtained. The particular set-up applies to well logging as well as to eddy current tube testing. In this work, a number of improvements are presented to existing models for an efficient calculation of the induced voltage. These include: domain truncation, novel treatment of arbitrary number of layers in order to avoid computational overflows and efficient time response calculation. The latter is based on a combination of Laplace inversion techniques for short and long time transient response.

# 1. Introduction

In transient or pulsed eddy-current testing (PECT), the excitation (transmitter) coil is driven by a pulsed current and the eddy current response from the inspected testpiece is sensed by another pickup (receiver) coil or by magnetic field sensors. Among other applications, PECT is used in monitoring pipe wall thinning since the growing field experience has demonstrated that PECT is effective for that matter [1], [2]. An important feature of the signal is its decay behavior which can be used for the quantitative evaluation of pipe wall thinning due to corrosion. The signal from a thinner pipe wall decays faster than that of a thicker one. Hence, there is a need for accurate evaluation of the response signal in the so called "long time domain" [3], [4].

The problem of theoretically analyzing PEC can be treated by first acknowledging that with respect to the geometry, testpiece and probe characteristics it is essentially a harmonic excitation problem. That is, the problem can be solved in the frequency domain and then the transient response can be derived either by Fourier superposition or by just replacing the $j\omega$ term with the Laplace variable $s$ and then apply Laplace inversion techniques. There is a tremendous bulk of analytical solutions for the harmonic excitation [5], therefore the focus should be on the effective inversion of the obtained Laplace expressions.

The theoretical study and analysis of PECT as applied to cylindrical structures such as insulated pipes is a matter of great interest. It is interesting that such studies are usually of numerical nature while analytical and semi-analytical ones limit themselves to the description of simplified configurations. The main simplification assumes a pipe diameter that is significantly larger than that of the excitation coil in order to model the configuration as a planar geometry [3]. When working with such a model it is easy to derive simple time constants from the signal decay rate in order to develop a method for thickness measurement [6].

Nevertheless, more representative axisymmetric models for the case of a multi-layer cylindrical system have been obtained in the frequency domain [7], followed by similar studies for pulsed excitations [8]-[9]. These use Fourier superposition, but such an approach can be time consuming and may suffer from the Gibb's phenomenon.

In this work, we also study an axisymmetric configuration of a transmitter-receiver coil system (otherwise known as reflection system) located in a multi-layer cylindrical conductor system. We then calculate the transient response of the voltage induced in the receiver coil by focusing on Laplace inversion [11] rather than on Fourier superposition.

Various methods can be used for the Laplace inversion depending on the characteristics of the expected time domain signal [12]. The use of Laplace transform requires an efficient method for inverting the frequency domain expressions to time domain. In the context of eddy current testing, a number of approaches have been used so far [13]-[16].

In addition to the Laplace treatment, we utilize domain truncation [17] and also present a novel systematic treatment of the multilayer cylindrical system that avoids possible issues with overflows owing to the modified Bessel functions behavior.

A critical examination of published studies is also presented in Sec.3 together with the description of the preferred methods for Laplace inversion.

## 2. Analysis
### 2.1 Induced voltage in the frequency domain

The eddy current interaction between a system of coaxial coils (transmitter and receiver) with an axisymmetric multilayer cylindrical conductor system, as shown in Fig.1, has been studied in [7]. The geometry extends infinitely in the $z$-direction and the harmonic excitation of the transmitter coil has the form $I\exp(j\omega t)$ where $I$ is the rms value of the excitation current and $\omega=2\pi f$ is the angular frequency. The induced voltage in the receiver coil, due to the cylindrical coil system, for arbitrary number of conductive layers $N$ is given by:

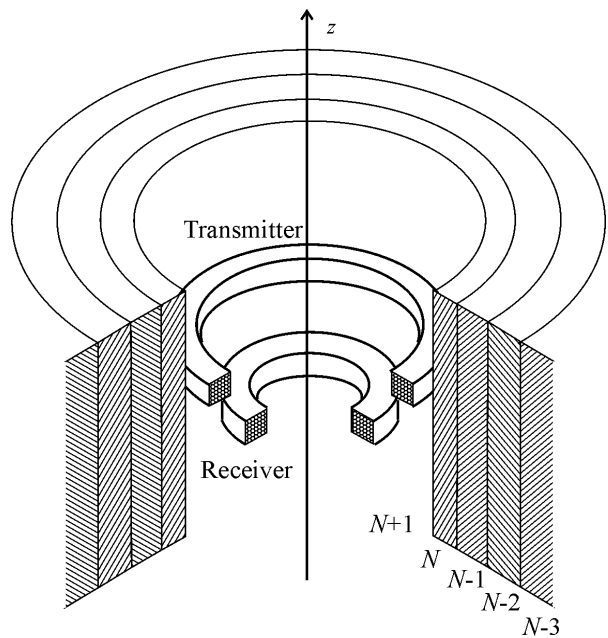

Fig. 1. Transmitter and receiver coil system in the presence of a multilayer cylindrical system.

$$V^R = j2\omega\mu_0 \frac{I_T w_T}{(r_{2T}-r_{1T})(z_{2T}-z_{1T})} \frac{w_R}{(r_{2R}-r_{1R})(z_{2R}-z_{1R})}$$

$$\times \int_0^\infty \frac{1}{q^6} \psi(qr_{1T},qr_{2T})\psi(qr_{1R},qr_{2R}) R(q) \qquad (1)$$

$$\times \left[\cos q(z_{1R}-z_{1T}) - \cos q(z_{2R}-z_{1T}) - \cos q(z_{1R}-z_{T}) + \cos q(z_{2R}-z_{2T})\right] dq$$

where

$$\psi(x_1,x_2) = \int_{x_1}^{x_2} x\, I_1(x)\, dx.$$

and $r_1$, $r_2$ are the coil radii, $z_1$, $z_2$ are the coil bottom and top height, $l = z_2 - z_1$ is the coil height, $w$ is the number of coil wire turns and the $T$ and $R$ subscripts refer to the Transmitter and Receiver coil respectively. Note also that the interface radii between the regions are denoted by $b_n$ and numbering of the cylindrical regions starts from the outer region and proceeds to the inner one. The air region where the coil lies is denoted as $N+1$. A layer designated as $n$ has conductivity $\sigma_n$ and relative magnetic permeability $\mu_{rn}$.

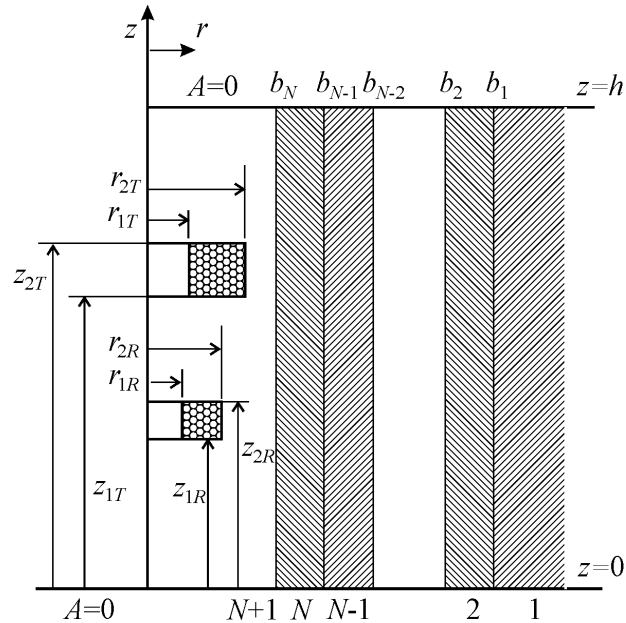

Fig. 2. Axisymmetric view of the truncated multicylindrical region with the two coils (transmitter and receiver) showing also region and interface numbering.

Another approach to the solution of the problem is the truncation of the solution domain from $z=0$ to $h$, as shown in Fig.2 and imposition of appropriate boundary conditions. In the examined configuration this involves Dirichlet conditions for the magnetic vector potential $A$, which in turn means magnetic insulation at the two boundaries. Following [17], the induced voltage in the receiver coil can be computed by an expression that involves a summation rather than the integration in (1)

$$V^R = \frac{j4\omega\mu_0\pi}{h} I^T \sum_{i=1}^{\infty} \frac{R(q_i)}{q_i^6} Y_T(q_i) Y_R(q_i) \tag{2}$$

where

$$Y(q_i) = \frac{w}{(r_2 - r_1)(z_2 - z_1)} \psi(q_i r_1, q_i r_2)\left[\cos(q_i z_1) - \cos(q_i z_2)\right] \tag{3}$$

with the discrete eigenvalue $q = \dfrac{i\pi}{h}$.

In both (1) and (2), the multilayer cylindrical system characteristics are described through the $R(q)$ term. These mathematical expressions are derived by applying the method of separation of variables in the Helmholtz equation for the magnetic vector potential. This potential is expressed in each cylindrical sub-region as a combination of sinusoidal and modified Bessel functions together with unknown coefficients. The latter are calculated by applying boundary conditions on the interfaces between the cylindrical sub-regions. Following [7], the $R(q)$ term can be calculated as follows:

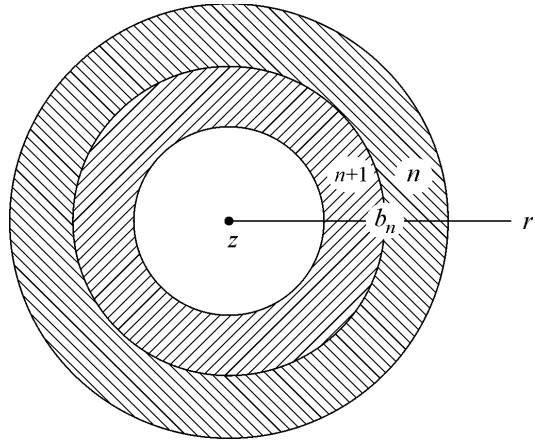

Fig. 3. Cross-section of the configuration for applying interface conditions between two cylindrical layers.

*2.2 The older approach*

For the arbitrary cylindrical layer $n$ the radial dependence of the magnetic vector potential is written as:

$$A_n(r) = C_n I_1(a_n r) + D_n K_1(a_n r) \tag{4}$$

where

$$a_n = \sqrt{q^2 + j\omega\mu_0\mu_{rn}\sigma_n}. \tag{5}$$

In the coil region, numbered $N+1$, the radial dependence of the magnetic vector potential is written as:

$$A_{N+1}(r) = C_{N+1}I_1(a_{N+1}r) + D_{N+1}K_1(a_{N+1}r) \tag{6}$$

where $D_{N+1}$ describes the effect of the coil in air and $C_{N+1}$ stands for the eddy current part from the cylindrical multilayer conductors. The relation between these terms can be written as $C_{N+1} = R(q)D_{N+1}$ where $R(q)$ is a reflection term, the same with the one showing up in (1). The procedure of calculating $R(q)$ starts by writing the coefficients of region $n+1$ in terms of the coefficients of region $n$ and for that matter, the interface conditions are applied on the boundary $b_n$. In terms of the radial dependence of the magnetic vector potential, these result in:

$$\begin{aligned} C_{n+1}I_1(a_{n+1}b_n) + D_{n+1}K_1(a_{n+1}b_n) &= C_nI_1(a_nb_n) + D_nK_1(a_nb_n) \\ C_{n+1}\beta_{n+1}I_0(a_{n+1}b_n) - D_{n+1}\beta_{n+1}K_0(a_{n+1}b_n) &= C_n\beta_nI_0(a_nb_n) - D_n\beta_nK_0(a_nb_n) \end{aligned} \tag{7}$$

and in matrix form they can be written as:

$$\begin{bmatrix} C_{n+1} \\ D_{n+1} \end{bmatrix} = \mathbf{T}_{n+1,n}\begin{bmatrix} C_n \\ D_n \end{bmatrix} \qquad \mathbf{T}_{n+1,n} = \frac{a_{n+1}b_n}{\beta_{n+1}}\begin{bmatrix} T_{11} & T_{12} \\ T_{21} & T_{22} \end{bmatrix} \tag{8}$$

where

$$\begin{aligned} T_{11} &= \beta_{n+1}K_0(a_{n+1}b_n)I_1(a_nb_n) + \beta_nI_0(a_nb_n)K_1(a_{n+1}b_n) \\ T_{12} &= \beta_{n+1}K_0(a_{n+1}b_n)K_1(a_nb_n) - \beta_nK_0(a_nb_n)K_1(a_{n+1}b_n) \\ T_{21} &= \beta_{n+1}I_0(a_{n+1}b_n)I_1(a_nb_n) - \beta_nI_0(a_nb_n)I_1(a_{n+1}b_n) \\ T_{22} &= \beta_{n+1}I_0(a_{n+1}b_n)K_1(a_nb_n) + \beta_nK_0(a_nb_n)I_1(a_{n+1}b_n) \end{aligned} \tag{9}$$

Starting from the outermost region ($n=1$) and proceeding towards the innermost ones

$$\begin{bmatrix} C_{N+1} \\ D_{N+1} \end{bmatrix} = \mathbf{T}_{N+1,N}\mathbf{T}_{N,N-1}\ldots\mathbf{T}_{n+1,n}\ldots\mathbf{T}_{3,2}\ldots\mathbf{T}_{2,1}\begin{bmatrix} C_1 \\ D_1 \end{bmatrix} \Rightarrow \begin{bmatrix} C_{N+1} \\ D_{N+1} \end{bmatrix} = \mathbf{U}\begin{bmatrix} C_1 \\ D_1 \end{bmatrix} \Rightarrow$$

$$\begin{bmatrix} C_{N+1} \\ D_{N+1} \end{bmatrix} = \begin{bmatrix} U_{11} & U_{12} \\ U_{21} & U_{22} \end{bmatrix}\begin{bmatrix} 0 \\ D_1 \end{bmatrix} \Rightarrow C_{N+1} = \frac{U_{12}}{U_{22}}D_{N+1} \Rightarrow R(q) = \frac{U_{12}}{U_{22}} \tag{10}$$

the final expression for the reflection term is derived.

## 2.3 The novel approach

In most cases, the analysis in Sec.2.2 works fine, but for a large number of layers or for extreme cases of material characteristics and/or frequency values may present overflow problems due to the behavior of Bessel functions for large arguments. With the approach described next, we avoid such problems. What we essentially do is to normalize each Bessel function with an argument at a radial distance $r$ in layer $n$ with the same Bessel function with an argument at either $b_n$ or $b_{n-1}$. For $I$-Bessel function we use $b_{n+1}$ while for $K$-Bessel functions we use $b_n$ so that the resulting ratio is always smaller than 1. Hence, in each region, the radial dependence of the magnetic vector potential is written as:

$$A_n(r) = C_n \frac{I_1(a_n r)}{I_1(a_n b_{n-1})} + D_n \frac{K_1(a_n r)}{K_1(a_n b_n)} \tag{11}$$

and especially for the inner region $N+1$ (where the coil belongs) $b_{N+1}=b_N$. For the region $n+1$, the radial dependence of the magnetic vector potential is written as:

$$A_{n+1}(r) = C_{n+1} \frac{I_1(a_{n+1} r)}{I_1(a_{n+1} b_n)} + D_{n+1} \frac{K_1(a_{n+1} r)}{K_1(a_{n+1} b_{n+1})} \tag{12}$$

Referring to Fig.3, imposition of interface conditions results in:

$$\begin{aligned}
& C_{n+1} + D_{n+1} \frac{K_1(a_{n+1} b_n)}{K_1(a_{n+1} b_{n+1})} = C_n \frac{I_1(a_n b_n)}{I_1(a_n b_{n-1})} + D_n \\
& C_{n+1}\beta_{n+1} \frac{I_0(a_{n+1} b_n)}{I_1(a_{n+1} b_n)} - D_{n+1}\beta_{n+1} \frac{K_0(a_{n+1} b_n)}{K_1(a_{n+1} b_{n+1})} = C_n\beta_n \frac{I_0(a_n b_n)}{I_1(a_n b_{n-1})} - D_n\beta_n \frac{K_0(a_n b_n)}{K_1(a_n b_n)} \Rightarrow \\
& C_{n+1}\beta_{n+1} \frac{I_0(a_{n+1} b_n)}{I_1(a_{n+1} b_n)} - D_{n+1}\beta_{n+1} \frac{K_0(a_{n+1} b_n)}{K_1(a_{n+1} b_n)} \frac{K_1(a_{n+1} b_n)}{K_1(a_{n+1} b_{n+1})} \\
& = C_n\beta_n \frac{I_0(a_n b_n)}{I_1(a_n b_n)} \frac{I_1(a_n b_n)}{I_1(a_n b_{n-1})} - D_n\beta_n \frac{K_0(a_n b_n)}{K_1(a_n b_n)}
\end{aligned} \tag{13}$$

If we then use the following abbreviations for the Bessel function ratios:

$$\begin{aligned}
I_n &= \frac{I_1(a_n b_n)}{I_1(a_n b_{n-1})} & K_n &= \frac{K_1(a_{n+1} b_n)}{K_1(a_{n+1} b_{n+1})} \\
I'_n &= \beta_n \frac{I_0(a_n b_n)}{I_1(a_n b_n)} & K'_n &= \beta_n \frac{K_0(a_n b_n)}{K_1(a_n b_n)} \\
I'_{n+1} &= \beta_{n+1} \frac{I_0(a_{n+1} b_n)}{I_1(a_{n+1} b_n)} & K'_{n+1} &= \beta_{n+1} \frac{K_0(a_{n+1} b_n)}{K_1(a_{n+1} b_n)}
\end{aligned} \tag{14}$$

the system of coefficients is written as:

$$\begin{aligned}
C_{n+1} + D_{n+1} K_n &= C_n I_n + D_n \\
C_{n+1} I'_{n+1} - D_{n+1} K'_{n+1} K_n &= C_n I'_n I_n - D_n K'_n
\end{aligned} \tag{15}$$

By solving (15) in terms of $C_{n+1}$ and $D_{n+1}$, the coefficients of region $n+1$ are written in terms of the coefficients of region $n$ as follows:

$$C_{n+1} = C_n \frac{(I'_n + K'_{n+1})I_n}{(I'_{n+1} + K'_{n+1})K_n} + D_n \frac{(K'_{n+1} - K'_n)}{(I'_{n+1} + K'_{n+1})}$$

$$D_{n+1} = C_n \frac{(I'_{n+1} - I'_n)I_n}{(I'_{n+1} + K'_{n+1})K_n} + D_n \frac{(I'_{n+1} + K'_n)}{(I'_{n+1} + K'_{n+1})K_n}$$

(16)

Now the new transfer matrix between the two regions is written as:

$$\begin{bmatrix} C_{n+1} \\ D_{n+1} \end{bmatrix} = \mathbf{T}_{n+1,n} \begin{bmatrix} C_n \\ D_n \end{bmatrix} \quad \mathbf{T}_{n+1,n} = \frac{1}{(I'_{n+1} + K'_{n+1})K_n} \begin{bmatrix} (I'_n + K'_{n+1})I_n & (K'_{n+1} - K'_n)K_n \\ (I'_{n+1} - I'_n)I_n & (I'_{n+1} + K'_n) \end{bmatrix}$$

(17)

and due to the normalization the final expression for the reflection term is modified to:

$$C_{N+1}I_1(a_{N+1}b_N) = \frac{U_{12}}{U_{22}} D_{N+1} K_1(a_{N+1}b_N) \Rightarrow C_{N+1} = \frac{U_{12}}{U_{22}} \frac{K_1(a_{N+1}b_N)}{I_1(a_{N+1}b_N)} D_{N+1} \Rightarrow R(q) = \frac{U_{12}}{U_{22}} \frac{K_1(a_{N+1}b_N)}{I_1(a_{N+1}b_N)}.$$

(18)

The procedure just described avoids all overflows of the expressions that may occur for a large number of cylindrical layers and for combinations of parameters that result in large Bessel arguments. This is possible because programming languages and mathematical packages provide computations for isolating the exponential behaviour of Bessel functions, for example in Matlab `besseli(v,z,1)` computes $I_v e^{-|\text{Re}\{z\}|}$ and `besselk(v,z,1)` computes $K_v e^z$ which means that the expressions $I_n, K_n$ in (14) can be computed by:

```
besseli(1,a_n*b_n,1)/ besseli(1,a_n*b_{n-1},1)*exp(real(a_n)*(b_n-b_{n-1}))
besselk(1,a_n*b_n,1)/ besselk(1,a_n*b_{n-1},1)*exp(-a_n*(b_n-b_{n-1}))
```

and the ratios in the expressions $I'_n, K'_n$ can be computed by:

```
besseli(0,a_n*b_n,1)/ besseli(1,a_n*b_n,1)
besselk(0,a_n*b_n,1)/ besselk(1,a_n*b_n,1).
```

Such methods for cancelling the exponential behaviour are especially useful in 3D configurations where ratios for multiple higher orders of Bessel functions are required.

## 3. Time response and Laplace inversion

As already stated, a universal method for the computation of the time response is Fourier superposition. However, too many frequencies are required for an accurate result. For example in [8] it is stated that a large number of 1500 distinct frequencies is required for a reliable computation. In order to reduce computational time, interpolation is used in the frequency spectrum and the number of frequencies is dictated by a logarithmic rule. In [9] such an interpolation is also utilized combined with an interpolation on the integration variable $q$ since the integral expression (1) is used rather than the summation one (2). The latter interpolation is justified in that paper due to the very large number of 1500 summation terms, referring to [8]. This is, however, a misconception since the number 1500 refers to the frequency values and not to the axial terms, thus making the use of tedious integral expression unnecessary.

In this paper, we focus on the inversion of the Laplace expressions. Having calculated $R(q_i)$, we return to (2) for the truncated domain, which describes also the induced EMF in the receiver coil, in the Laplace domain, when the excitation current is pulsed. This is done by replacing $j\omega$ with the Laplace variable $s$.

$$\tilde{V}^R(s) = \frac{4\mu_0\pi}{h}\tilde{I}^T(s)\sum_{i=1}^{\infty}\frac{s\tilde{R}(q_i,s)}{q_i^6}Y_T(q_i)Y_R(q_i) \qquad (19)$$

For a step current $I(t) = I_0 u(t)$, the Laplace transform is $\tilde{I}(s) = I_0/s$, so the problem now is the computation of the inverse Laplace transform of $\tilde{R}(q_i,s)$, with discrete eigenvalues $q_i$. The inverse Laplace of (19) is written as

$$V^R(t) = \frac{4\mu_0\pi}{h}I_0\sum_{i=1}^{\infty}\frac{R(q_i,t)}{q_i^6}Y_T(q_i)Y_R(q_i) \qquad (20)$$

where $R(q_i,t)$ is the inverse Laplace transform of $\tilde{R}(q_i,s)$.

*3.1 The numerical method*

A universal approach can be based on the numerical inversion of the Laplace transform (NILT), methods for which rank among ones which are widely used for time-domain simulations. From many developed methods, those based on FFT and ε-algorithm for the acceleration of series convergence seem to be convenient from the point of view of both desired speed and accuracy [18]. The publicly available Matlab code in [18] has been used for fast computation. The algorithm requires only a number of sampling points and the time interval in which the transient signal is to be computed (starting from *t*=0). The specific code is very efficient since it utilizes the advantageous feature of Matlab language to run in parallel on multidimensional arrays without necessity to use outer loop structures. This means that time response calculations can be performed simultaneously for all eigenvalues $q_i$ and leads to essential savings in CPU time.

*3.2 The Stefhest method*

Apart from universal numerical approaches, there are many other algorithms to invert the Laplace transform, each of them having specific advantages and disadvantages depending on the form of the function to be inverted. For the short time transient signal that has the form of exponential damping (as in our case) we can utilize the Stehfest algorithm which involves specific time instances. The time response is computed by a sum that involves weighted values of the Laplace function:

$$f(t) \approx \frac{\ln 2}{t} \sum_{i=1}^{n} V_i F\left(\frac{\ln 2}{t} i\right) \qquad (21)$$

where the weights $V_i$ are determined by:

$$V_i = (-1)^{\frac{n}{2}+1} \sum_{k=\frac{i+1}{2}}^{\min\left(i,\frac{n}{2}\right)} \frac{k^{\frac{n}{2}}(2k)!}{\left(\frac{n}{2}-k\right)!k!(k-1)!(i-k)!(2k-i)!} \qquad (22)$$

An optimal choice for the number of terms $n$ in the summation of (21) is $10 \leq n \leq 14$, here we are using 14. Note that (22) is incorrectly written in [10], which can attributed to the fact that the same expression is also incorrectly written in [12] from which it was taken. The Stehfest algorithm gives very good results in the short time domain. In Sec.4 we will study its reliable use by comparing its results with the universal numerical method.

*3.3 Pole extraction*

The residue theorem is an alternative method for Laplace inversion and involves the computation of poles of the function to be inverted. The method is ideal for the so called "long time" calculation since in this case only a few poles need to be found for the accurate description of the transient signal. The diffusive nature of the eddy current phenomenon implies that all poles are real, leading thus to distinct decay modes for the acquired signal. Suitable methods are utilized for bracketing the poles prior to their exact calculation, following in part the method developed in [15]. A similar approach has been adopted also in [14] for a coil encircling a conductive rod.

The problem now is the computation of the inverse Laplace transform of the reflection term $\tilde{R}(q_i, s)$ which can be written as:

$$\tilde{R}(q_i, s) = \frac{R_1(q_i, s)}{R_2(q_i, s)} \qquad (23)$$

For each eigenvalue $q_i$ we seek the poles $s_k$ of (23), then by invoking the Heaviside expansion theorem the time dependent term is written as:

$$R(q_i, t) = \sum_k \frac{R_1(q_i, s_k)}{R'_2(q_i, s_k)} e^{s_k t} \qquad (24)$$

Since the poles $s_k$ have to be real negative a search routine that starts from $s=0$ and proceeding to the real negative axis can be easily implemented. Such routines are found in every mathematical package.

An important feature of the pole extraction method is that only a limited number of poles needs to be used. Especially the long time response may be computed with just one pole! This feature is studied in [15] and is also utilized in [19] in order to describe the behavior of the eddy current system by just one pole that corresponds to the dominant mode of the decay rate of the transient signal. The feature could also be used in [16] instead of using a large constant number of poles.

## 4. Results and Discussion

The examined system comprises two concentric tubes, each one 10mm thick, with an air gap between them. This creates a multilayer cylindrical system with $N=4$ layers where layers numbered 1 and 3 are air. The tubes are assumed to be of the same material and are either carbon steel (ferromagnetic) or stainless steel (non ferromagnetic). The cylindrical system and coil data are given in Table I:

Table I. Coils and tubes characteristics

|  | Transmitter | Receiver | Layers radii | | Tubes material: steel | |
|---|---|---|---|---|---|---|
| $r_1$ [mm] | 20 | 20 | $b_1$ [mm] | 70 | $\sigma$ [MS/m] | 3 |
| $r_2$ [mm] | 30 | 30 | $b_2$ [mm] | 60 | $\mu_r$ | 100 (carbon) |
| $l$ [mm] | 40 | 10 | $b_3$ [mm] | 50 | | |
| $wt$ | 1600 | 10000 | $b_4$ [mm] | 40 | $\sigma$ [MS/m] | 3 |
| $L_{air}$ [mH] | 83.24 | 5738.5 | | | $\mu_r$ | 1 (stainless) |

Given that the axial gap between the coils is set to $g=10$mm, the required coil distances in Fig.2 are calculated from:

$$z_{1T} = \frac{h}{2}, \quad z_{2T} = z_{1T} + l_T, \quad z_{2R} = \frac{h}{2} - g, \quad z_{1R} = z_{2R} - l_R$$

Results are provided for the induced voltage in the receiver coil when the transmitter coil is excited by a step current of 1A. The transient responses are shown in Fig.4 for the three inversion methods described in Sec.3. The logarithmic scale in the voltage axis depicts better the long time behaviour of the transient signal. Adopting as reference the NILT method we can easily observe that the pole extraction behaves very well in the long time domain while the Stehfest method behaves well for the short time domain and starts to deviate from a certain point in time.

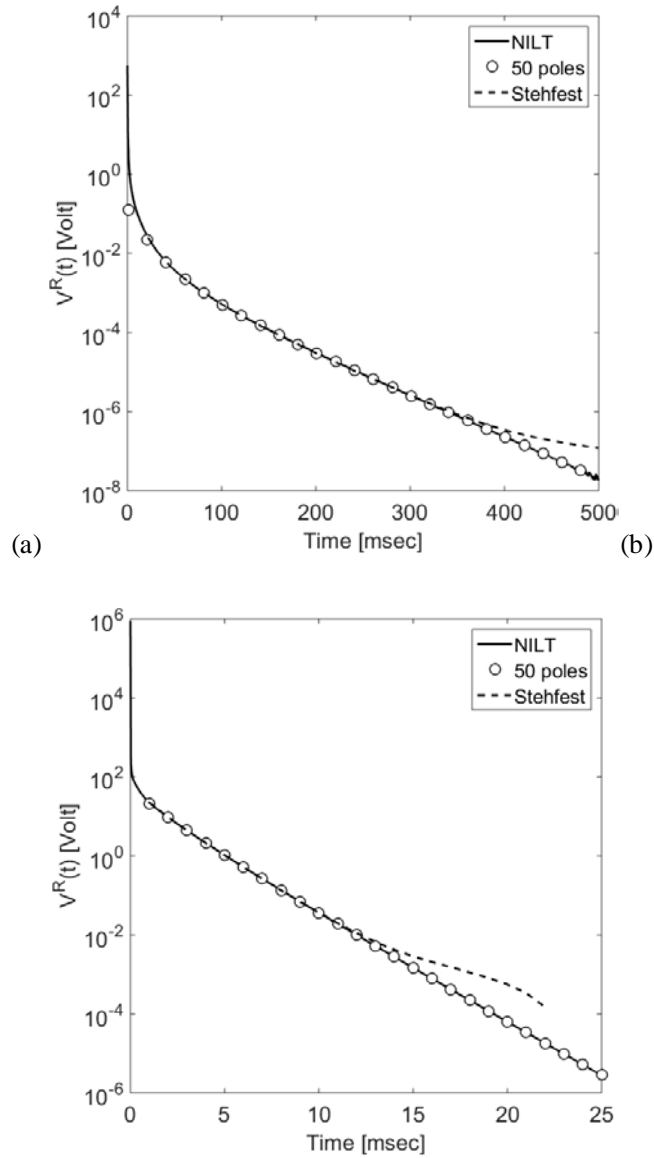

Fig. 4. Long time comparison of Stehfest and pole extraction methods to NILT for (a) Carbon and (b) Stainless Steel. Only the first pole is used for each $q_i$ eigenvalue.

It seems then that these two methods, Stehfest and pole extraction, are complementary and can be used in the short and long time domain, respectively. The point of transition that we have established empirically, after many simulations, for a system of tubes of the same material is the time defined by $t_m = \mu_0 \mu_r \sigma b^2$, where $b$ stands for the total thickness of the tubes.

From Fig.4, it is evident that the logarithmic signal curve is steep in the case of the stainless steel due to its smaller relative magnetic permeability and hence to the weaker eddy currents induced in the tube wall. This is exactly the expected behaviour [3], [6], [19].

In the summation of (2) we have used 50 terms, and for the Laplace inversion we only use 1 pole per eigenvalue (hence 50 poles). The boundary of the truncated domain is set to $h=100 r_{2T}$ for carbon and $h=20 r_{2T}$ for stainless steel. The fact that we use only 1 pole for each eigenvalue is a crude computation that nevertheless gives good results for the long time. For the long time domain, the number of eigenvalues and thus the number of poles can be decreased further. We have run a parametric study regarding the number of poles that can be used to reliably compute the induced voltage. Fig.5 shows results that compare pole extraction method (with decreasing number of poles) to the exact NILT solution. It is clear that the limiting case of using just 1 pole can also be used for a reliable representation in the long time region.

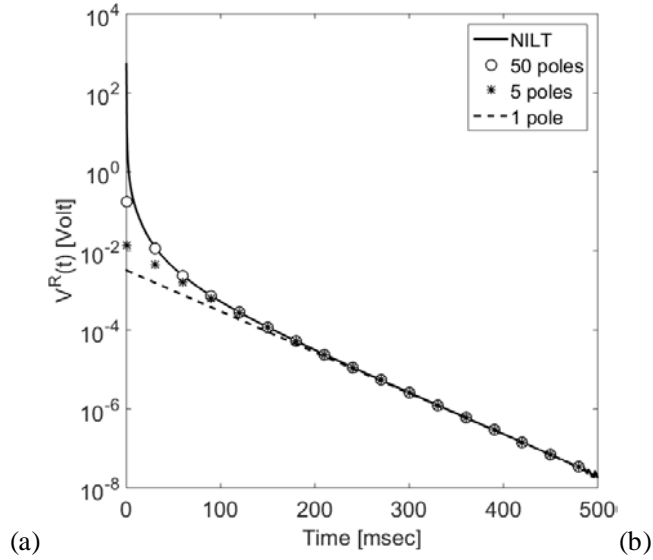

(a)

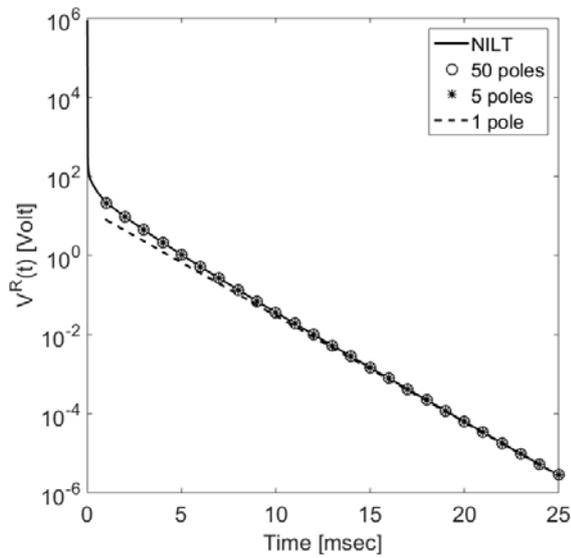

(b)

Fig. 5. Long time comparison of NILT to pole extraction method with different number of poles for (a) Carbon and (b) Stainless Steel.

In this case, the induced voltage is computed by only the first term of (20), which for the cases of the carbon (CS) and stainless steel (SS) tubes and the data in Table reduces to:

$$V_{CS}^R(t) = 0.0032 e^{-23.87t} \Rightarrow \log_{10} V_{CS}^R(t) = -2.49 - 10.4t$$
$$V_{SS}^R(t) = 14.717 e^{-619.5t} \Rightarrow \log_{10} V_{SS}^R(t) = 1.17 - 269t$$
(25)

These essentially define the dominant decay modes of the induced voltage signal.

Fig. 6 shows the signals for wall thinning of 50% on each side of the two carbon steel tubes together with the signal in the case of the absence of a tube. Obviously the decay rate (curve inclination) changes with the amount of wall thinning and its specific location, i.e. in which tube it happens and if it happens in the ID or OD of the specific tube. Not all cases can be resolved in the long time domain. For example, the case of the absence of one tube cannot be easily distinguished (which tube is absent) since the total remaining thickness is the same. In this case we refer to the short time behaviour, shown in Fig.6b where the voltage is plotted in linear scale. These transient signals are computed with the Stehfest method.

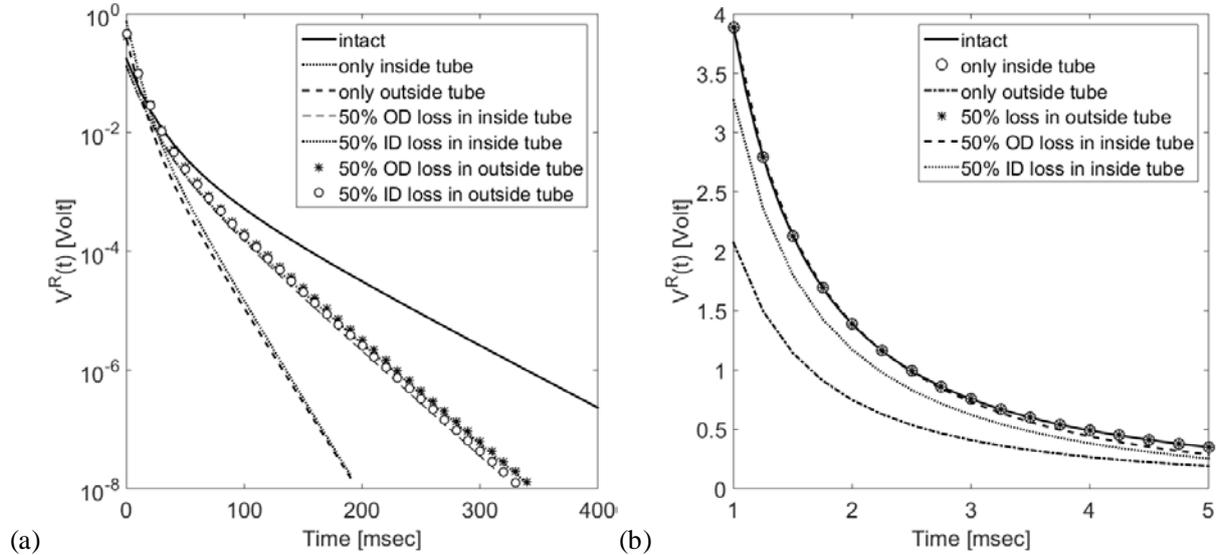

Fig. 6. Signals of wall thickness loss and for each tube for carbon steel for (a) long time and (b) short time.

Hence, the thickness and location of the wall loss can be derived from the combination of amplitude in the short time domain and decay rate (curve inclination) in the long time domain.

From the depicted results it is easily derived that the pole extraction can be used for calibration of instruments that utilize PEC for thickness measurements. The fact that the presented model describes the exact cylindrical configuration improves the method accuracy since no planar geometry simplification is utilized.

Nevertheless, as shown in experiments by Ulapane et al in [4] and [6], if we consider the typical PEC waveform the noise "destroys" the response when the signal is a few orders of magnitude smaller than the initial value. From this perspective, high accuracy for long time may be useless especially when the transition point (from using Stehfest to using pole extraction) lies in the time area affected by noise. In this case, only the use of Stehfest method, which is better for short time response, is enough.

## 4. Conclusions

We have presented a number of improvements on existing models for the transient response of a coil system inside a cylindrical conductive configuration. The truncation of the solution domain offers certain computational advantages since it simplifies the induced voltage calculation. Normalization by using Bessel function ratios enables an efficient solution when the number of cylindrical layers increases, it may also improve the computations for a 3D configuration when higher order Bessel functions are required. Finally it was shown that regarding the transient signal computation: For the short time, the Stehfest method is very efficient while for the long time the residue theorem can give a quick and exact computation with a minimum number of poles. In any case numerical methods based on Fourier method together combined with series acceleration algorithm provides a universal and reliable alternative. For long time computations it was shown that even one (1) pole is adequate for the whole computation since it represents the dominant decay mode and can be used for tube wall thickness evaluations.

Possible extensions of the presented improvements include the 2D axisymmetric geometry with encircling coil around a multilayer cylindrical conductor system. Furthermore, the 3D geometries can be treated in a similar manner since relevant solutions exist in the frequency domain for both ID and OD coils as well as for eccentrically placed tubes.

**References**

[1] A. Sophian, G.Y. Tian and M. Fan, Pulsed eddy current nondestructive testing and evaluation: A review, Chin. J. Mech. Eng. 30 (2017), 500-514.

[2] Q. Luo, Y. Shi, Z. Wang, W. Zhang and Y. Li, A study of applying pulsed remote field eddy current in ferromagnetic pipes testing, Sensors 17 (2017), 1-8.

[3] W. Cheng, Pulsed eddy current testing of carbon steel pipes' wall-thinning through insulation and cladding, Journal of NDE 31 (2012), 215-224.

[4] N. Ulapane and L.Nguyen, Review of pulsed eddy current signal feature extraction methods for conductive ferromagnetic material-thickness quantification, Electronics 8 (2019), 470.

[5] Theodoulidis T.P., Kriezis E.E., Eddy current canonical problems (with applications to nondestructive evaluation), Tech Science Press (2006).

[6] N. Ulapane, A. Alempijevic, J.V. Miro and T. Vidal-Valleja, Nondestructive evaluation of ferromagnetic material thickness using Pulsed Eddy Current sensor detector coil voltage decay rate, NDT and E International, 100 (2018), 108-114.

[7] Dodd C.V., Cheng C.C., Deeds W.E., "Induction coils coaxial with an arbitrary number of cylindrical conductors", Journal of Applied Physics, Vol.45, No.2, pp.638-647, 1974.

[8] Y. Li, X. Liu, Z. Chen, H. Zhao and W. Cai, A fast forward model of pulsed eddy current inspection of multilayered tubular structures, International Journal of Applied Electromagnetics and Mechanics 45 (2014), 417-423.

[9] Z. Xue, M. Fan, B. Cao and D. Wen, A fast numerical method for the analytical model of pulsed eddy current for pipelines, Insight 62 (2020), 1-7.

[10] M. Fan, B. Cao and Y. Wang, Computation of coil-induced voltage due to a defect-free plate using Stehfest's method for pulsed eddy current evaluation, Insight 52 (2010), 1-3.

[11] A.M. Cohen, Numerical methods for Laplace transform inversion, Springer (2007).

[12] H. Hassanzadeh and M. Pooladi-Darvish, Comparison of different Laplace inversion methods for engineering applications, Applied Mathematics and Computation 189 (2007), 1966-1981.

[13] F. Fu and J.R. Bowler, Transient eddy current driver pickup probe response due to a conductive plate, IEEE Transactions on Magnetics 42 (2006), 2029-2037.

[14] F. Fu and J.R. Bowler, Transient eddy current response due to a conductive cylindrical rod, AIP Conference Proceedings 894 (2007), 332.

[15] T. Theodoulidis, Developments in calculating the transient eddy-current response from a conductive plate, IEEE Transactions on Magnetics 44 (2008), 1894–1896.

[16] X. Chen and Y. Lei, Time-domain analytical solutions to pulsed eddy current field excited by a probe coil outside a conducting ferromagnetic pipe, NDT&E International 68 (2014), 22-27.

[17] Bowler J.R., Theodoulidis T.P., "Eddy currents induced in a conducting rod of finite length by a coaxial encircling coil", Journal of Physics D: Applied Physics, Vol.38, pp.2861-2868, 2005.

[18] L. Brancik, Programs for fast numerical inversion of Laplace transforms in Matlab language environment, http://dsp.vscht.cz/konference_matlab/MATLAB99/.

[19] X. Chen and Y. Lei, Excitation current waveform for eddy current testing on the thickness of ferromagnetic plates, NDT&E International 66 (2014), 28-33.